\documentclass[reqno,12pt,draft]{amsart}
\usepackage{eurosym}
\usepackage{amsmath,amsthm,amsfonts,amssymb}
\usepackage[english]{babel}

\newtheorem{theorem}{Theorem}
\newtheorem{lemma}{Lemma}
\newtheorem{remark}{Remark}

\newtheorem{proposition}{Proposition}
\input{tcilatex}

\begin{document}
\date{}

\begin{center}
{\LARGE Boundary control and homogenization: optimal climatization through
smart double skin boundaries }
\end{center}

\vspace{0.3cm}

\begin{center}
{\large J.I. D\'{\i}az }$^{1}${\large , A.V. Podolskiy }$^{2}${\large \ and
T.A. Shaposhnikova }$^{2}${\large \ }

($^{1}$) Instituto de Matematica Interdisciplinar (IMI)

Dpto. An\'{a}lisis Matem\'{a}tico y Matem\'{a}tica Aplicada, Universidad
Complutense de Madrid, 28040 Madrid, Spain

jidiaz@ucm.es,

($^{2}$) Lomonosov Moscow State University

Department of Differential Equations

Faculty of Mechanics and Mathematics, Leninskie Gory, 119991, GSP-1, Moscow,
Russian Federation

originalea@ya.ru \ \ 

shaposh.tan@mail.ru

\bigskip

\textit{Dedicated to the unforgettable Olga A. Oleinik (1925-2001)}
\end{center}

\vspace{0.3cm}

\textbf{Abstract}. We consider the homogenization of an optimal control
problem in which the control $v$ is placed on a part $\Gamma _{0}$ of \ the
boundary and the spatial domain contains a thin layer of \textquotedblleft
small particles\textquotedblright , very close to the controlling boundary,
and a Robin boundary condition is assumed on the boundary of those
\textquotedblleft small particles\textquotedblright . This problem can be
associated with the climatization modeling of \textit{Bioclimatic Double
Skin Fa\c{c}ades} which was developed in modern architecture as a tool for
energy optimization. We assume that the size of the particles and the
parameters involved in the Robin boundary condition are critical (and so
they justify the occurrence of some \textquotedblleft strange
terms\textquotedblright\ in the homogenized problem). The cost functional is
given by a weighted balance of the distance (in a $H^{1}$-type metric) to a
prescribed target internal temperature $u_{T}$ and the proper cost of the
control $v$ (given by its $L^{2}(\Gamma _{0})$ norm). We prove the (weak)
convergence of states ${u_{\varepsilon }}$ and of the controls $%
v_{\varepsilon }$ to some functions, ${u_{0}}$ and $v_{0}$, respectively,
which are completely identified: ${u_{0}}$ satisfies an artificial boundary
condition on $\Gamma _{0}$ and $v_{0}$ is the optimal control associated to
a limit cost functional $J_{0}$ in which the \textquotedblleft boundary
strange term\textquotedblright\ on $\Gamma _{0}$ arises. This information on
the limit problem makes much more manageable the study of the optimal
climatization of such double skin structures.

\textrm{\textit{Key words: \textrm{optimal control, }homogenization, thin
layer of \textquotedblleft small particles\textquotedblright , critical
case, cost functional convergence }}

AMS Subject Classification: 35B27, 93C20, 49N05

\textrm{\textit{\vspace{0.3cm} }}

\section{\textrm{\textit{\textbf{Introduction}}}}

A well-known energy optimization technique in modern architecture is the
theory of\textit{\ smart fa\c{c}ade systems (}also called as\textit{\
Bioclimatic Double Skin Fa\c{c}ades) }in which\textit{\ }climatization takes
place by means of active glass windows (see, e.g., \cite{Braham}, \cite{Ben
Bonham} and \cite{Azad}). Today, there are different types of active glass
in the market: LCD Liquid Crystal, Gasochromic, SPO suspended particles,
Electrochromic, etc. See, e.g., the case of fluids and windows in \cite{Ben
Bonham} and \cite{Claros-Padial}.

From the mathematical view point, many different climatization models have
been proposed in the literature: see, for instance, Chapter 1 of the
excellent book by Duvaut and Lions \cite{Duvaut-Lions}. Some studies on
internal climatization and homogenization can be found in \cite{Timofte}. In
this paper, we will analyze a simplified formulation of \textit{Double Skin
Fa\c{c}ades }in which there is an active flux control $v_{\varepsilon }$,
located in a part $\Gamma _{0}$ of the boundary, and a kind of \textit{%
celosia (}latticed windows called in this way in Spanish) traditionally made
of masonry, wood, or a combination of these materials. We assume that the 
\textit{celos\'{\i}a} is formed by a set of periodical small thermostats, of
period $\varepsilon >0$, located in an internal thin layer located very
close to the controlling boundary. So, $\varepsilon $ represents a small
parameter related to the characteristic \textit{celosia}.

Our simplified optimal control problem assumes that the state of the system $%
{u_{\varepsilon }}$ (the internal temperature) satisfies a Poisson equation
in the internal domain $\Omega _{\varepsilon }$ of $\mathbb{R}^{n}$, with $%
n\geq 3$, which is defined as the external domain to the set of periodical
small thermostats (here represented by a set of $\varepsilon $-periodically
balls) on whose contours $S_{\varepsilon }$ a given climatization law
(represented by a Robin boundary condition with a large parameter $%
\varepsilon ^{-\gamma }$ as coefficient, where $\gamma =\frac{n-1}{n-2}$)
takes place. Non-symmetrical shapes, and/or the case $n=2$, can also be
considered thanks to the techniques presented in \cite{DPSH Dokladi 5}, but
for the sake of simplicity in the presentation we will not develop it here.
We assume that any thermostat has a critical radius $a_{\varepsilon }$,
where $a_{\varepsilon }=C_{0}\varepsilon ^{\alpha }$ and $\alpha =\frac{n-1}{%
n-2}$. As in many other frameworks (see many examples and references in the
monograph \cite{DGCSh 1 Book}), this critical size leads to the occurrence
of strange terms in the homogenized problem (in contrast with what happens
for other possible sizes).

It is assumed that the cost functional $J_{\varepsilon }(v_{\varepsilon })$
is given by a weighted balance of the distance (in a $H^{1}(\Omega
_{\varepsilon })$ type metric) to a prescribed target internal temperature $%
u_{T}$ and the proper cost of the control $v$ (given by its $L^{2}(\Gamma
_{0})$ norm). Our main result proves the (weak) convergence of solutions ${%
u_{\varepsilon }}$ and of the controls $v_{\varepsilon }$ to some functions, 
${u_{0}}$ and $v_{0},$ respectively, which are completely identified: ${u_{0}%
}$ satisfies the Poisson equation in the whole domain $\Omega $ (which we
assume to be a bounded open set with $\partial \Omega $ of class $C^{1}$)
and $v_{0}$ is the boundary optimal control but in an \textit{artificial }%
boundary condition in which the thermostats effects are located on the own
controllability boundary $\Gamma _{0}.$ Moreover, we prove the convergence
of the cost functional $J_{\varepsilon }(v_{\varepsilon })$ to a new cost
functional $J_{0}(v_{0})$ in which the boundary \textit{strange term }arises%
\textit{. }This information on the limit problem makes much more manageable
the study of the optimal climatization of such double skin structures.

In the last section, we consider the pure homogenization process (without
any control, $v_{\varepsilon }\equiv 0$) and prove the convergence of the
energies: an information which is stronger than the mere weak convergence ${%
u_{\varepsilon }\rightharpoonup u_{0}}$ in $H^{1}(\Omega )$.

This paper complements the scope considered by previous papers in the
literature concerning optimal control problems in which the controls are
located in different parts of the spatial domain (see \cite{Saint Pau-Zou 2}%
, \cite{Stromqvuist 10}, \cite{Pod-Sh 2020 3}, \cite{DPSH 2021 4} and \cite%
{DPSH Dokladi 5}).

\section{\textrm{\textit{\textbf{Problem statement. Adjoint optimality
problem }}}}

Let $\Omega \subset \mathbb{R}_{+}^{n}=\{x\in \mathbb{R}^{n}:x_{n}>0\}$ be a
bounded open set of $\mathbb{R}^{n}$, $n\geq 3$, with $\partial \Omega $ of
class $C^{2}$, and $\partial \Omega =\Gamma _{0}\cup \Gamma _{1}$, is
assumed to be of class $C^{1}$, where $\Gamma _{0}=\partial \Omega \cap
\{x_{n}=0\}\neq \emptyset $ is the $(n-1)$-dimensional domain on the plane $%
x_{n}=0$ which represents the controlling boundary, $\Gamma _{1}=\partial
\Omega \setminus \overline{\Gamma _{0}}$. Define $Y_{0}=(-1/2,1/2)^{n-1}%
\times (0,1)$, $Y_{\varepsilon }^{j}=\varepsilon Y_{0}+\varepsilon j$, $%
j=(j_{1},\ldots j_{n-1},0)$, $j_{i}\in \mathbb{Z},i=1,\ldots n-1$. We denote
by $P_{\varepsilon }^{j}$ the center of the cube $Y_{\varepsilon }^{j}$, $%
G_{\varepsilon }^{j}=a_{\varepsilon }G_{0}+\varepsilon j$, where $G_{0}$ is
the unit ball with the center coinciding with the center $(0,\ldots 0,1/2)$
of the cube $Y_{0}$ and $a_{\varepsilon }=C_{0}\varepsilon ^{\alpha }$ with $%
\alpha =(n-1)/(n-2)$. We define $\Upsilon _{\varepsilon }=\{j\in \mathbb{Z}%
^{n}:j=(j_{1},\ldots ,j_{n-1},0),Y_{\varepsilon }^{j}\subset \Omega \}.$ It
is easy to see (as in Chapter 6 of \cite{DGCSh 1 Book}) that $|\Upsilon
_{\varepsilon }|\cong d\varepsilon ^{1-n}$, $d=const>0$, where $|\Upsilon
_{\varepsilon }|$ denotes the cardinality of the set of isolated points $%
\Upsilon _{\varepsilon }$.

We introduce the sets 
\begin{equation*}
G_{\varepsilon }=\bigcup_{j\in \Upsilon _{\varepsilon }}G_{\varepsilon
}^{j},\,\,S_{\varepsilon }=\bigcup_{j\in \Upsilon _{\varepsilon }}\partial {%
G_{\varepsilon }^{j}},\,\,\Omega _{\varepsilon }=\Omega \setminus \overline{%
G_{\varepsilon }},\,\,\partial \Omega _{\varepsilon }=S_{\varepsilon
}\bigcup \Gamma _{0}\bigcup \Gamma _{1}.\newline
\end{equation*}%
The set $G_{\varepsilon }$ represents the \textit{celos\'{\i}a} or double
skin. It is localized as a subset of $\Omega \cap \{x\in \mathbb{R}%
^{n}:x_{n}\in (0,\varepsilon )\}$ (see Figure 1).

\FRAME{ftbpFU}{4.7409in}{3.9807in}{0pt}{\Qcb{Example of a spatial domain
with a\textit{\ double skin boundary.}}}{}{celosia3.jpg}{\special{language
"Scientific Word";type "GRAPHIC";maintain-aspect-ratio TRUE;display
"USEDEF";valid_file "F";width 4.7409in;height 3.9807in;depth
0pt;original-width 8.5314in;original-height 7.1563in;cropleft "0";croptop
"1";cropright "1";cropbottom "0";filename '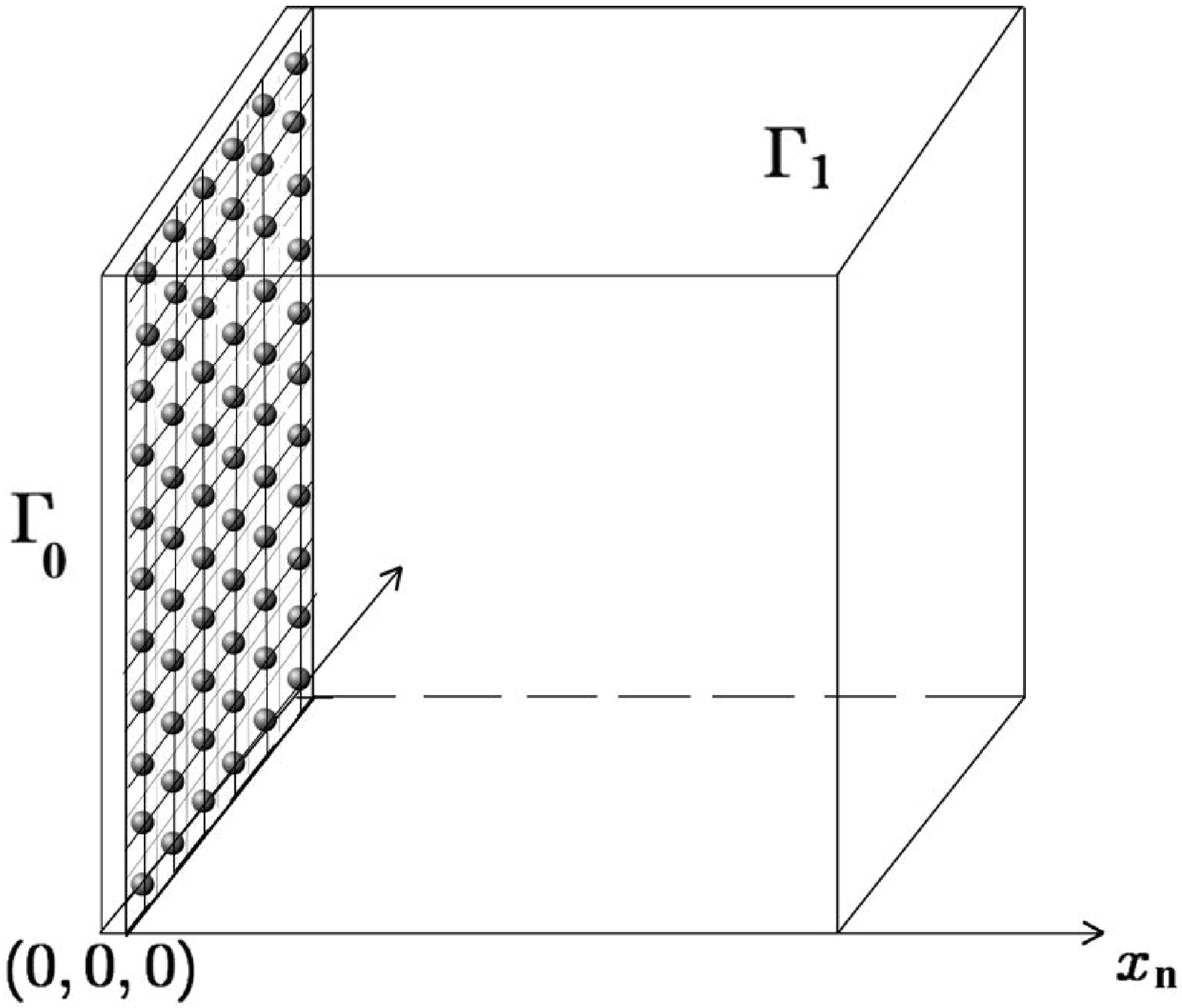';file-properties
"XNPEU";}}

For an arbitrary function $v\in L^{2}(\Gamma _{0})$, we denote by $%
u_{\varepsilon }(v)\in H^{1}(\Omega _{\varepsilon },\Gamma _{1})$ the
solution of the problem

\begin{equation}
\left\{ 
\begin{array}{lr}
-\Delta {u_{\varepsilon }(v)}=f, & x\in \Omega _{\varepsilon }, \\ 
\partial _{\nu }u_{\varepsilon }(v)+{\varepsilon }^{-\gamma
}a(x)u_{\varepsilon }(v)=0, & x\in S_{\varepsilon }, \\ 
\partial _{\nu }u_{\varepsilon }(v)=v, & x\in \Gamma _{0}, \\ 
u_{\varepsilon }(v)=0, & x\in \Gamma _{1},%
\end{array}%
\right.  \label{Eq 1}
\end{equation}%
where $f\in L^{2}(\Omega )$, $a(x)\in C^{\infty }(\overline{\Omega })$, $%
a(x)\geq a_{0}=const>0$, and the notation $\partial _{\nu }g$ represents the
partial derivative along the outward unit normal vector $\nu $ to the
boundary. Here, the space $H^{1}(\Omega _{\varepsilon },\Gamma _{1})=\left\{
w\in H^{1}(\Omega _{\varepsilon })\text{ such that }w=0\text{ on }\Gamma
_{1}\right\} .$

We assume to be given a target function $u_{T}\in H^{1}(\Omega )$ and we
consider the cost functional

\begin{equation*}
J_{\varepsilon }:L^{2}(\Gamma _{0})\rightarrow \mathbb{R},\newline
\end{equation*}%
given by 
\begin{equation}
J_{\varepsilon }(v)=\frac{\eta }{2}\int_{\Omega _{\varepsilon }}{B(x)}\nabla
(u_{\varepsilon }(v)-u_{T})\nabla (u_{\varepsilon }(v)-u_{T})dx+\frac{N}{2}%
\Vert {v}\Vert _{L^{2}(\Gamma _{0})}^{2},\newline
\label{cost functional eps 3}
\end{equation}%
where the weighted balance is defined through the arbitrary positive
constants $\eta ,N$, and the $H^{1}-$metric is defined by the $n\times n$
symmetric matrix $B(x)=(b_{ij}(x))$ such that 
\begin{equation}
\lambda _{1}|\xi |^{2}\leq b_{ij}(x)\xi _{i}\xi _{j}\leq \lambda _{2}|\xi
|^{2},  \label{coercivity B 2}
\end{equation}%
for a.e. $x\in \Omega $, $\lambda _{i}=const>0$, $i=1,2$, $B\in C^{1}(%
\overline{\Omega })^{n\times n}$.

By well-known results (see e.g. \cite{Lions book 9}, \cite{Fursikov}, \cite%
{Trolstz}), there exists a unique optimal control $v_{\varepsilon }\in
L^{2}(\Gamma _{0})$

\textrm{\textit{%
\begin{equation}
J_{\varepsilon }(v_{\varepsilon })=\min\limits_{v\in L^{2}(\Gamma
_{0})}J_{\varepsilon }(v).  \label{Optim control eps 4}
\end{equation}%
}}One of the main goals of this paper is to find the limit as $\varepsilon
\rightarrow 0$ of the optimal control $v_{\varepsilon }$, of the associate
state $u_{\varepsilon }(v_{\varepsilon })$ and of the cost functional value $%
J_{\varepsilon }(v_{\varepsilon })$.

In order to characterize the optimal control $v_{\varepsilon }$, we will
study the particularization of the abstract version of the Pontryagin
maximum principle applied to elliptic PDEs mentioned in Section 1.3 of Lions 
\cite{Lions book 9}.

\begin{proposition}
\label{Propo1}\textit{Let }$u_{T}\in H^{1}(\Omega ),$ $v_{\varepsilon }\in
L^{2}(\Gamma _{0})$\textit{\ and }$u_{\varepsilon }(v_{\varepsilon })\in
H^{1}(\Omega _{\varepsilon },\Gamma _{1})$\textit{\ be the target state, the
optimal control and the associate optimal state, respectively. Let }$%
P_{\varepsilon }\in $\textit{\ }$H^{1}(\Omega _{\varepsilon },\Gamma _{1})$ 
\textit{be the unique solution of the problem}%
\begin{equation}
\left\{ 
\begin{array}{lr}
\Delta {P_{\varepsilon }}=div(B(x)\nabla (u_{\varepsilon }-u_{T}), & x\in
\Omega _{\varepsilon }, \\ 
\partial _{\nu }P_{\varepsilon }-(B(x)\nabla (u_{\varepsilon }-u_{T}),\nu )+{%
\varepsilon }^{-\gamma }a(x)P_{\varepsilon }=0, & x\in S_{\varepsilon }, \\ 
\partial _{\nu }P_{\varepsilon }-(B(x)\nabla (u_{\varepsilon }-u_{T}),\nu
)=0, & x\in \Gamma _{0}, \\ 
P_{\varepsilon }=0, & x\in \Gamma _{1}.%
\end{array}%
\right.  \label{probl P eps 5}
\end{equation}%
\textit{Then, the optimal control }$v_{\varepsilon }$ \textit{is given by}%
\begin{equation}
v_{\varepsilon }=-\frac{\eta }{N}P_{\varepsilon }\mathit{.}
\label{Optimal control  and p eps}
\end{equation}%
\textit{\ }
\end{proposition}

\noindent \textit{Proof}. Since $v_{\varepsilon }$ is the optimal control,
we know that for any other control $v\in L^{2}(\Gamma _{0})$ 
\begin{equation*}
\lim\limits_{\lambda \rightarrow 0}\frac{1}{\lambda }(J_{\varepsilon
}(v_{\varepsilon }+\lambda {v})-J_{\varepsilon }(v_{\varepsilon }))=0.
\end{equation*}%
It is easy to see that if, for a given $\lambda \in \mathbb{R}$, we define%
\begin{equation*}
\theta _{\varepsilon }=\frac{1}{\lambda }(u_{\varepsilon }(v_{\varepsilon
}+\lambda {v})-u_{\varepsilon }(v_{\varepsilon })),
\end{equation*}%
then $\theta _{\varepsilon }$ is a weak solution to the problem

\textrm{\textit{%
\begin{equation}
\left\{ 
\begin{array}{lr}
\Delta {\theta _{\varepsilon }}=0, & x\in \Omega _{\varepsilon }, \\ 
\partial _{\nu }\theta _{\varepsilon }+{\varepsilon }^{-\gamma }a(x)\theta
_{\varepsilon }=0, & x\in S_{\varepsilon }, \\ 
\partial _{\nu }\theta _{\varepsilon }=v, & x\in \Gamma _{0}, \\ 
\theta _{\varepsilon }=0, & x\in \Gamma _{1}.%
\end{array}%
\right.  \label{Probl theta eps 6}
\end{equation}%
}}Then, we have

\begin{equation}
0=\lim\limits_{\lambda \rightarrow 0}\frac{1}{\lambda }\Bigl(J_{\varepsilon
}(v_{\varepsilon }+\lambda {v})-J_{\varepsilon }(v_{\varepsilon })\Bigr)%
=\eta \int_{\Omega _{\varepsilon }}B(x)\nabla (u_{\varepsilon
}(v_{\varepsilon })-u_{T})\nabla {\theta _{\varepsilon }}dx+N\int_{\Gamma
_{0}}v_{\varepsilon }{v}d\hat{x}.\newline
\label{limit eps 7}
\end{equation}%
We recall that if $P_{\varepsilon }\in H^{1}(\Omega _{\varepsilon },\Gamma
_{1})$ is a solution of (\ref{probl P eps 5}), then we have the integral
identity 
\begin{equation}
\int_{\Omega _{\varepsilon }}\nabla {P_{\varepsilon }}\nabla \phi {dx}+{%
\varepsilon }^{-\gamma }\int_{S_{\varepsilon }}a(x)P_{\varepsilon }\phi {ds}%
=\int_{\Omega _{\varepsilon }}B(x)\nabla (u_{\varepsilon }-u_{T})\nabla \phi 
{dx},\newline
\label{Integral iden 8}
\end{equation}%
for an arbitrary function $\phi \in H^{1}(\Omega _{\varepsilon },\Gamma
_{1}) $. Then, we can set $\phi =\theta _{\varepsilon }$ in it and use $%
P_{\varepsilon }$ as a test function in the integral identity for the
problem (\ref{Probl theta eps 6}). Subtracting the one from the other, we
get 
\begin{equation*}
\eta \int_{\Omega _{\varepsilon }}B(x)\nabla (u_{\varepsilon
}(v_{\varepsilon })-u_{T})\nabla \theta _{\varepsilon }dx+N\int_{\Gamma
_{0}}v_{\varepsilon }v{d\hat{x}}=\int_{\Gamma _{0}}(\eta {P_{\varepsilon }}+N%
{v_{\varepsilon }})v{d\hat{x}}=0,
\end{equation*}%
from which we deduce that $v_{\varepsilon }=-\dfrac{\eta }{N}P_{\varepsilon }
$.$_{\blacksquare }$

\bigskip

In order to get the homogenization (as $\varepsilon \rightarrow 0$) we will
use the usual $H^{1}$-extensions of functions $u_{\varepsilon }$ and $%
P_{\varepsilon }$ to $\Omega $ which we denote by $\widetilde{u_{\varepsilon
}}$ and $\widetilde{P_{\varepsilon }}$ (see, e.g., Section 3.1.1 of \cite%
{DGCSh 1 Book} and its references).

Then, from the properties of the extension operator (see~\cite{Zu-Sh 2011 7}%
) and estimates (\ref{estim direct 16}), we have

\begin{theorem}
\label{Theorem 1}\textit{\ }Let $f\in L^{2}(\Omega ),$ $u_{T}\in
H^{1}(\Omega )$ and let $(u_{\varepsilon },P_{\varepsilon })\in H^{1}(\Omega
_{\varepsilon },\Gamma _{1})^{2}$ be the weak solution of the coupled system%
\textrm{\textit{%
\begin{equation}
\left\{ 
\begin{array}{lr}
-\Delta {u_{\varepsilon }}=f, & x\in \Omega _{\varepsilon }, \\ 
\Delta {P_{\varepsilon }}=div(B(x)\nabla (u_{\varepsilon }-u_{T})), & x\in
\Omega _{\varepsilon }, \\ 
\partial _{\nu }u_{\varepsilon }+{\varepsilon }^{-\gamma }a(x)u_{\varepsilon
}=0, & x\in S_{\varepsilon }, \\ 
\partial _{\nu }P_{\varepsilon }-(B(x)\nabla (u_{\varepsilon }-u_{T}),\nu )+{%
\varepsilon }^{-\gamma }a(x)P_{\varepsilon }=0, & x\in S_{\varepsilon }, \\ 
\partial _{\nu }u_{\varepsilon }=-\dfrac{\eta }{N}P_{\varepsilon }, & x\in
\Gamma _{0}, \\ 
\partial _{\nu }P_{\varepsilon }-(B(x)\nabla (u_{\varepsilon }-u_{T}),\nu
)=0, & x\in \Gamma _{0}, \\ 
u_{\varepsilon }=P_{\varepsilon }=0, & x\in \Gamma _{1}.%
\end{array}%
\right.  \label{Coupled system eps 9}
\end{equation}%
}}

\noindent Then, 
\begin{equation}
\Vert \widetilde{u_{\varepsilon }}\Vert _{H^{1}(\Omega ,\Gamma _{1})}\leq
C,\,\,\Vert \widetilde{P_{\varepsilon }}\Vert _{H^{1}(\Omega ,\Gamma
_{1})}\leq C,\newline
\label{Apriori estim 17}
\end{equation}%
\noindent and thus there exists some subsequences (still denoted as the
original ones) such that 
\begin{equation}
\widetilde{u_{\varepsilon }}\rightharpoonup u_{0},\,\,\widetilde{%
P_{\varepsilon }}\rightharpoonup P_{0},\,\,\,\mbox{weakly in}%
\,\,H^{1}(\Omega ,\Gamma _{1}),  \label{Weak conv 18}
\end{equation}

\noindent as $\varepsilon \rightarrow 0,$ for some $(u_{0},P_{0})\in $ $%
H^{1}(\Omega ,\Gamma _{1})^{2}.$
\end{theorem}

\noindent \textit{Proof}. We start by getting some a priori estimates for $%
u_{\varepsilon }$ and $P_{\varepsilon }$. From the integral identity for the
function $P_{\varepsilon }$, we derive 
\begin{equation}
\begin{array}{c}
\displaystyle \int_{\Omega _{\varepsilon }}|\nabla {P_{\varepsilon }}|^{2}dx+%
{\varepsilon }^{-\gamma }\int_{S_{\varepsilon }}a(x)P_{\varepsilon
}^{2}ds=\int_{\Omega _{\varepsilon }}B(x)\nabla (u_{\varepsilon
}-u_{T})\nabla {P_{\varepsilon }}dx\leq \\[.3cm] 
\leq \dfrac{1}{2}\Vert \nabla {P_{\varepsilon }}\Vert _{L^{2}(\Omega
_{\varepsilon })}^{2}+C\Vert \nabla (u_{\varepsilon }-u_{T})\Vert
_{L^{2}(\Omega _{\varepsilon })}^{2}.\newline
\end{array}
\label{Inequ 10}
\end{equation}%
The constant $C$ doesn't depend on $\varepsilon $ here and below. From here,
we conclude 
\begin{equation}
\Vert \nabla {P_{\varepsilon }}\Vert _{L^{2}(\Omega _{\varepsilon })}^{2}+{%
\varepsilon }^{-\gamma }\Vert {P_{\varepsilon }}\Vert _{L^{2}(S_{\varepsilon
})}^{2}\leq C\Vert \nabla (u_{\varepsilon }-u_{T})\Vert _{L^{2}(\Omega
_{\varepsilon })}^{2}.  \label{sec estim 11}
\end{equation}%
From the integral identity for the function $u_{\varepsilon }$, we have 
\begin{equation}
\int_{\Omega _{\varepsilon }}\nabla {u_{\varepsilon }}\nabla {P_{\varepsilon
}}dx+{\varepsilon }^{-\gamma }\int _{S_{\varepsilon }}a(x)u_{\varepsilon
}P_{\varepsilon }ds=\int_{\Omega _{\varepsilon }}f{P_{\varepsilon }}dx-\frac{%
\eta }{N}\int_{\Gamma _{0}}P_{\varepsilon }^{2}d\hat{x}.  \label{Equal 12}
\end{equation}%
From the integral identity for the function $P_{\varepsilon }$, we get 
\begin{equation}
\int_{\Omega _{\varepsilon }}\nabla {u_{\varepsilon }}\nabla {P_{\varepsilon
}}dx+{\varepsilon }^{-\gamma }\int_{S_{\varepsilon }}a(x)P_{\varepsilon
}u_{\varepsilon }ds=\int_{\Omega _{\varepsilon }}B(x)\nabla (u_{\varepsilon
}-u_{T})\nabla {u_{\varepsilon }}dx.  \label{Sec Equa 13}
\end{equation}%
Subtracting the equality (\ref{Equal 12}) from (\ref{Sec Equa 13}), we get 
\begin{equation}
\frac{\eta }{N}\int_{\Gamma _{0}}P_{\varepsilon }^{2}d\hat{x}+\int_{\Omega
_{\varepsilon }}B(x)\nabla (u_{\varepsilon }-u_{T})\nabla (u_{\varepsilon
}-u_{T})dx=\int_{\Omega _{\varepsilon }}f{P_{\varepsilon }}dx-\int_{\Omega
_{\varepsilon }}B(x)\nabla {u_{T}}\nabla (u_{\varepsilon }-u_{T})dx.
\label{third equ 14}
\end{equation}%
This equality and the condition (\ref{coercivity B 2}) imply

\begin{equation}
\Vert {P_{\varepsilon }}\Vert _{L^{2}(\Gamma _{0}))}\leq C(\Vert \nabla {%
u_{T}}\Vert _{L^{2}(\Omega )}+\Vert {f}\Vert _{L^{2}(\Omega )})\text{, }%
\Vert \nabla (u_{\varepsilon }-u_{T})\Vert _{L^{2}(\Omega _{\varepsilon
})}\leq C(\Vert \nabla {u_{T}}\Vert _{L^{2}(\Omega )}+\Vert {f}\Vert
_{L^{2}(\Omega )}).\newline
\label{estima 15}
\end{equation}%
From (\ref{sec estim 11}) and (\ref{estima 15}), we derive 
\begin{equation}
\Vert {u_{\varepsilon }}\Vert _{H^{1}(\Omega _{\varepsilon },\Gamma
_{1})}\leq C,\,\,\Vert {P_{\varepsilon }}\Vert _{H^{1}(\Omega _{\varepsilon
},\Gamma _{1})}\leq C.\newline
\label{estim direct 16}
\end{equation}%
Since $\widetilde{u_{\varepsilon }}$ and $\widetilde{P_{\varepsilon }}$, are
the $H^{1}$-extensions of $u_{\varepsilon }$ and $P_{\varepsilon }$ to $%
\Omega $, from the properties of the extension operator (see, e.g., ~\cite%
{Zu-Sh 2011 7}) and estimates (\ref{estim direct 16}), we get the estimates (%
\ref{Apriori estim 17}) and thus we conclude the weak convergences indicated
in (\ref{Weak conv 18}) for some subsequences. $_{\blacksquare }$

\bigskip

In order to identify the limit problem satisfied by the pair $(u_{0},P_{0})$%
, we need several auxiliary results.

\section{Auxiliary statements}

For $j\in \Upsilon _{\varepsilon }$, we introduce the function $%
w_{\varepsilon }^{j}(x)$ as being the unique~solution to the capacity
boundary value problem 
\begin{equation}
\left\{ 
\begin{array}{lr}
\Delta {w_{\varepsilon }^{j}}=0, & x\in T_{{\varepsilon }/4}^{j}\setminus 
\overline{G_{\varepsilon }^{j}}, \\ 
w_{\varepsilon }^{j}=1, & x\in \partial {G_{\varepsilon }^{j}}, \\ 
w_{\varepsilon }^{j}=0, & x\in \partial {T_{{\varepsilon }/4}^{j},}%
\end{array}%
\right.  \label{Auxiliary probl 19}
\end{equation}%
where $T_{{\varepsilon }/4}^{j}$ is the ball of radius ${\varepsilon }/4$
with center in the point $P_{\varepsilon }^{j}$ (the center of the cube $%
Y_{\varepsilon }^{j}$). It is easy to see (Section 3.1.5.1 of \cite{DGCSh 1
Book}) that 
\begin{equation*}
w_{\varepsilon }^{j}(x)=\frac{|x-P_{\varepsilon }^{j}|^{2-n}-\Bigl(\frac{%
\varepsilon }{4}\Bigr)^{2-n}}{a_{\varepsilon }^{2-n}-\Bigl(\frac{\varepsilon 
}{4}\Bigr)^{2-n}}.\newline
\end{equation*}%
Define the extension function 
\begin{equation}
W_{\varepsilon }(x)=\left\{ 
\begin{array}{ll}
w_{\varepsilon }^{j}(x), & x\in T_{{\varepsilon }/4}^{j}\setminus \overline{%
G_{\varepsilon }^{j}},\,j\in \Upsilon _{\varepsilon }, \\ 
1, & x\in \bigcup _{j\in \Upsilon _{\varepsilon }}G_{\varepsilon }^{j}, \\%
[.2cm] 
0, & x\in \Omega \setminus \overline{\bigcup _{j\in \Upsilon _{\varepsilon }}%
}T_{{\varepsilon }/4}^{j}.%
\end{array}%
\right.  \label{Extension W ep 20}
\end{equation}%
It is clear that we have $W_{\varepsilon }\rightharpoonup 0,$ weakly in $%
H_{0}^{1}(\Omega ),$ as $\varepsilon \rightarrow 0$.

\begin{lemma}
Let $\phi \in C^{1}(\overline{\Omega },\Gamma _{1})$. Then,

\begin{equation}
\lim_{{\varepsilon }\rightarrow 0}\int_{\Omega _{\varepsilon }}(B\nabla {%
W_{\varepsilon }},\nabla {W_{\varepsilon }})\phi {dx}=\frac{%
C_{0}^{n-2}(n-2)\omega _{n}}{n}\int_{\Gamma _{0}}tr{B(0,\hat{x})}\phi (0,%
\hat{x}){d\hat{x}},\newline
\label{extensi Test func  21}
\end{equation}%
\noindent where $trB(x)=\sum\limits_{j=1}^{n}b_{jj}(x)$ is the trace of the
matrix $B(x)$ and $\omega _{n}$ is the surface area of the unit sphere in $%
\mathbb{R}^{n}$.
\end{lemma}

\noindent \textit{Proof.} Note that if $\widetilde{B}$ is the matrix with
the constant elements (with respect to the $y-$variable), then

\begin{equation}
\int_{S_{1}^{0}}(\widetilde{B}y,y)ds=\int_{T_{1}^{0}}div(\widetilde{B}y)dy=tr%
{\widetilde{B}}|T_{1}^{0}|=tr{\widetilde{B}}\frac{\omega _{n}}{n},\newline
\label{matrix B 22}
\end{equation}%
where $S_{1}^{0}$, $T_{1}^{0}$ are the unit sphere and the unit ball with
the center at the coordinate's origin.

\noindent Using equality (\ref{matrix B 22}), if $r=\sqrt{%
\sum\limits_{j=1}^{n}(x_{i}-P_{\varepsilon ,i}^{j})^{2}},$ we get 
\begin{equation*}
\begin{array}{l}
\displaystyle \lim_{\varepsilon \rightarrow 0}\int_{\Omega _{\varepsilon
}}(B(x)\nabla {W_{\varepsilon }},\nabla {W_{\varepsilon }})\phi
(x)dx=\lim_{\varepsilon \rightarrow 0}\sum_{j\in \Upsilon _{\varepsilon
}}\phi (P_{\varepsilon }^{j})\int_{T_{{\varepsilon }/4}^{j}\setminus 
\overline{G_{\varepsilon }^{j}}}(B(P_{\varepsilon }^{j})\nabla {%
w_{\varepsilon }^{j}},\nabla {w_{\varepsilon }^{j}})dx= \\[.5cm] 
\displaystyle =\lim_{{\varepsilon }\rightarrow 0}\sum_{j\in \Upsilon
_{\varepsilon }}a_{\varepsilon }^{2(n-2)}(n-2)^{2}\phi (P_{\varepsilon
}^{j})\sum_{k,l=1}^{n}\int_{T_{{\varepsilon }/4}^{j}\setminus \overline{%
G_{\varepsilon }^{j}}}b_{kl}(P_{\varepsilon }^{j})(x_{l}-P_{\varepsilon
,l}^{j},x_{k}-P_{{\varepsilon },k}^{j}){r^{-2n}}dr= \\[.5cm] 
\displaystyle =\lim_{\varepsilon \rightarrow 0}\sum_{j\in \Upsilon
_{\varepsilon }}a_{\varepsilon }^{2(n-2)}(n-2)^{2}\phi (P_{\varepsilon
}^{j})\int_{a_{\varepsilon }}^{{\varepsilon }/4}r^{1-n}dr%
\int_{S_{1}^{0}}(B(P_{\varepsilon }^{j})y,y)ds= \\[.5cm] 
\displaystyle =\lim_{\varepsilon \rightarrow 0}\sum _{j\in \Upsilon
_{\varepsilon }}a_{\varepsilon }^{2(n-2)}(n-2)^{2}\frac{\omega _{n}}{n}\phi
(P_{\varepsilon }^{j})trB(P_{\varepsilon }^{j})\int_{a_{\varepsilon }}^{{%
\varepsilon }/4}r^{1-n}dr= \\ 
\displaystyle =C_{0}^{n-2}(n-2)\frac{\omega _{n}}{n}\lim_{{\varepsilon }%
\rightarrow 0}\sum _{j\in \Upsilon _{\varepsilon }}\phi (P_{\varepsilon
}^{j})trB(P_{\varepsilon }^{j}){\varepsilon }^{n-1}=\frac{C_{0}^{n-2}\omega
_{n}(n-2)}{n}\int_{\Gamma _{0}}trB(0,\hat{x})\phi (0,\hat{x})d\hat{x},%
\end{array}%
\end{equation*}%
which proves the result.$_{\blacksquare }$

\bigskip

We consider now the auxiliary boundary value problem

\textrm{\textit{%
\begin{equation}
\left\{ 
\begin{array}{lr}
\Delta {h_{\varepsilon }}=div(B\nabla {W_{\varepsilon }}), & x\in \Omega
_{\varepsilon }, \\ 
\partial _{\nu }h_{\varepsilon }-(B\nabla {W_{\varepsilon }},\nu )+{%
\varepsilon }^{-\gamma }a(x)h_{\varepsilon }=0, & x\in S_{\varepsilon }, \\ 
\partial _{\nu }h_{\varepsilon }-(B\nabla {W_{\varepsilon }},\nu )=0, & x\in
\Gamma _{0}, \\ 
h_{\varepsilon }=0, & x\in \Gamma _{1}.%
\end{array}%
\right.  \label{Problem h 23}
\end{equation}%
}}As usual, we say that function $h_{\varepsilon }\in H^{1}(\Omega
_{\varepsilon },\Gamma _{1})$ is a weak solution of problem (\ref{Problem h
23}) if it satisfies the integral identity

\begin{equation}
\int _{\Omega _{\varepsilon }}\nabla {h_{\varepsilon }}\nabla \phi {dx}+{%
\varepsilon }^{-\gamma }\int_{S_{\varepsilon }}a(x)h_{\varepsilon }\phi {ds}%
=\int_{\Omega _{\varepsilon }}(B\nabla {W_{\varepsilon }},\nabla \phi )dx,
\label{Integral iden h 24}
\end{equation}%
where $\phi $ is an arbitrary function from the space $H^{1}(\Omega
_{\varepsilon },\Gamma _{1})$. Setting $\phi =h_{\varepsilon }$ in (\ref%
{Auxiliary probl 19}) and using the properties of function $W_{\varepsilon }$%
, we derive an estimate $\Vert \widetilde{h}_{\varepsilon }\Vert
_{H^{1}(\Omega ,\Gamma _{1})}\leq C$. Therefore, there exists a subsequence
(we preserve the notation of the original one) such that 
\begin{equation}
\widetilde{h}_{\varepsilon }\rightharpoonup h_{0},\,\,\mbox{weakly in}%
\,\,H^{1}(\Omega ,\Gamma _{1}),\,\,{\varepsilon }\rightarrow 0.
\label{weak conver h 25}
\end{equation}

The following theorem identifies the limit function $h_{0}.$

\begin{theorem}
The function $h_{0}$ defined in (\ref{weak conver h 25}) is a weak solution
of the boundary value problem 
\begin{equation}
\left\{ 
\begin{array}{lr}
\Delta {h_{0}}=0, & x\in \Omega , \\ 
\partial _{\nu }h_{0}+\mathcal{A}_{1}\dfrac{a(x)}{a(x)+C_{n}}h_{0}=-\mathcal{%
A}_{2}\dfrac{trB(x)a(x)}{a(x)+C_{n}}, & x\in \Gamma _{0}, \\ 
h_{0}=0, & x\in \Gamma _{1},%
\end{array}%
\right.  \label{Limi proble h0 26}
\end{equation}%
\noindent where $\mathcal{A}_{1}(n)=(n-2)C_{0}^{n-2}\omega _{n}$, $\mathcal{A%
}_{2}(n)=\dfrac{\mathcal{A}_{1}(n)}{n}$, $C_{n}=\dfrac{n-2}{C_{0}}$.
\end{theorem}

\noindent \textit{Proof. }We take as test function in the integral identity (%
\ref{Integral iden h 24}) the function $\phi =W_{\varepsilon }\psi $, where $%
\psi \in C^{\infty }(\overline{\Omega },\Gamma _{1})$ (this is called an 
\textit{oscillating test function }according to the general Tartar's method)
and we get

\begin{equation}
\int_{\Omega _{\varepsilon }}\nabla {W_{\varepsilon }}\nabla (h_{\varepsilon
}\psi )dx+{\varepsilon }^{-\gamma }\int_{S_{\varepsilon }}a(x)h_{\varepsilon
}\psi {ds}=\int_{\Omega _{\varepsilon }}(B(x)\nabla {W_{\varepsilon }}%
,\nabla {W_{\varepsilon }})\psi {dx}+\beta _{\varepsilon },\newline
\label{Oscilant 27}
\end{equation}%
where $\beta _{\varepsilon }\rightarrow 0$, $\varepsilon \rightarrow 0$.

\noindent From (\ref{Oscilant 27}), we derive%
\begin{equation}
\begin{array}{c}
\displaystyle \sum_{j\in \Upsilon _{\varepsilon }}\int_{\partial {%
G_{\varepsilon }^{j}}}\partial _{\nu }w_{\varepsilon }^{j}h_{\varepsilon
}\psi {ds}+{\varepsilon }^{-\gamma }\int_{S_{\varepsilon
}}a(x)h_{\varepsilon }\psi {ds}+ \\[.6cm] 
\displaystyle +\sum\limits_{j\in \Upsilon _{\varepsilon }}\int _{\partial {%
T_{{\varepsilon }/4}^{j}}}\partial _{\nu }w_{\varepsilon }^{j}h_{\varepsilon
}\psi {ds}=\frac{C_{0}^{n-2}(n-2)\omega _{n}}{n}\int_{\Gamma _{0}}trB(\hat{x}%
)\psi (\hat{x})d\hat{x}.%
\end{array}
\label{oscilant bis 28}
\end{equation}

\noindent Hence, it follows that

\begin{equation}
\begin{array}{c}
\displaystyle \varepsilon ^{-\gamma }\int_{S_{\varepsilon }}\left (a(x)+%
\frac{n-2}{C_{0}}\right )h_{\varepsilon }\psi {ds}= \\[.5cm] 
\displaystyle-\sum _{j\in \Upsilon _{\varepsilon }}\int_{\partial {T_{{%
\varepsilon }/4}^{j}}}\partial _{\nu }w_{\varepsilon }^{j}h_{\varepsilon
}\psi {ds}+\frac{C_{0}^{n-2}\omega _{n}(n-2)}{n}\int_{\Gamma _{0}}trB(\hat{x}%
)\psi (\hat{x})d\hat{x}+\beta _{\varepsilon }.%
\end{array}
\label{oscil tris 29}
\end{equation}%
We set $\psi (x)=\dfrac{a(x)}{a(x)+C_{n}}v(x)$, where $v\in C^{\infty }(%
\overline{\Omega },\Gamma _{1})$, $C_{n}=\dfrac{n-2}{C_{0}}$, and get 
\begin{equation}
\begin{array}{c}
\displaystyle \lim _{\varepsilon \rightarrow 0}{\varepsilon }^{-\gamma
}\int_{S_{\varepsilon }}a(x)h_{\varepsilon }v{ds}= \\[.4cm] 
\displaystyle =C_{0}^{n-2}(n-2)\omega _{n}\int_{\Gamma _{0}}\frac{a(\hat{x})%
}{a(\hat{x})+C_{n}}h_{0}{v}d\hat{x}+\frac{C_{0}^{n-2}(n-2)\omega _{n}}{n}%
\int_{\Gamma _{0}}\frac{trB(\hat{x})a(\hat{x})}{a(\hat{x})+C_{n}}{v}(\hat{x}%
)d\hat{x}.\newline
\end{array}
\label{Oscilant4 30}
\end{equation}

\noindent Therefore, from (\ref{Integral iden h 24}) and (\ref{Oscilant4 30}%
), we obtain the integral identity for the function $h_{0}$

\begin{equation*}
\displaystyle \int _{\Omega }\nabla {h_{0}}\nabla \phi {dx}+\mathcal{A}%
_{1}\int _{\Gamma _{0}}\frac{a(\hat{x})}{a(\hat{x})+C_{n}}h_{0}\phi {d\hat{x}%
}+ \displaystyle +\frac{\mathcal{A}_{1}}{n}\int _{\Gamma _{0}}\frac{trB(\hat{%
x})a(\hat{x})}{a(\hat{x})+C_{n}}\phi (\hat{x})d\hat{x}=0.
\end{equation*}

\noindent This implies the statement of the theorem. $_{\blacksquare }$

\section{The main result}

So, in order to obtain the characterization of $u_{0}$ and $P_{0}$, we have
to pass to the limit in the identity (\ref{Integral iden 8}).

\begin{theorem}
\label{Theorem Main}Let $n\geq 3$, $\alpha =\gamma =\frac{n-1}{n-2}$ and let 
$(u_{\varepsilon },P_{\varepsilon })$ be a weak solution of the coupled
system (\ref{Coupled system eps 9}). Then, the pair $(u_{0},P_{0})$ defined
in (\ref{Weak conv 18}) is a weak solution of the system 
\begin{equation}
\left\{ 
\begin{array}{lr}
-\Delta {u_{0}}=f, & x\in \Omega , \\ 
\Delta {P_{0}}=div(B\nabla (u_{0}-u_{T})), & x\in \Omega , \\ 
\partial _{\nu }u_{0}+\mathcal{A}_{1}\dfrac{a(x)}{a(x)+C_{n}}u_{0}=-\dfrac{%
\eta }{N}P_{0}, & x\in \Gamma _{0}, \\ 
\partial _{\nu }P_{0}-(B(x)\nabla (u_{0}-u_{T}),\nu )+\mathcal{A}_{1}\dfrac{%
a(x)}{a(x)+C_{n}}P_{0}-\mathcal{A}_{2}\dfrac{trB(x)a^{2}(x)}{(a(x)+C_{n})^{2}%
}u_{0}=0, & x\in \Gamma _{0}, \\ 
u_{0}=P_{0}=0, & x\in \Gamma _{1},%
\end{array}%
\right.  \label{Coupled limit system 31}
\end{equation}%
where $trB(x)=\sum\limits_{j=1}^{n}b_{jj}(x)$ is the trace of the matrix $%
B(x)$, $\mathcal{A}_{1}(n)=(n-2)C_{0}^{n-2}\omega _{n}$, $\mathcal{A}_{2}(n)=%
\dfrac{\mathcal{A}_{1}(n)}{n}$, $C_{n}=\dfrac{n-2}{C_{0}}$ and $\omega _{n}$
is the surface area of the unit sphere in $\mathbb{R}^{n}$. In addition, if
we define the functional 
\begin{equation*}
\begin{array}{ll}
J_{0}(v) & \displaystyle \hspace*{-.2cm}=\frac{\eta }{2}\int_{\Omega
}B\nabla (u_{0}(v)-u_{T})\nabla (u_{0}(v)-u_{T})dx+\frac{N}{2}\int_{\Gamma
_{0}}v^{2}d\hat{x} \\[.5cm] 
& \displaystyle +\frac{\eta {\mathcal{A}_{1}}}{2n}\int_{\Gamma _{0}}trB(\hat{%
x})\Bigl(\frac{a(\hat{x})}{a(\hat{x})+C_{n}}\Bigr)^{2}u_{0}^{2}(v)d\hat{x},%
\end{array}%
\end{equation*}
then%
\begin{equation}
\lim _{{\varepsilon }\rightarrow 0}J_{\varepsilon }(v_{\varepsilon
})=J_{0}(v_{0}),  \label{definition J0 41}
\end{equation}
where $v_{0}=-\dfrac{\eta }{N}P_{0}$ on $\Gamma _{0}$. In particular, $%
v_{\varepsilon }\rightharpoonup $ $v_{0}$ weakly in $L^{2}(\Gamma _{0})$,
and $v_{0}$ is the optimal control of the problem\textrm{\textit{%
\begin{equation}
J_{0}(v_{0})=\inf _{v\in L^{2}(\Gamma _{0})}J_{0}(v),\newline
\label{Limit control problem}
\end{equation}%
}}associated to the state problem

\begin{equation*}
\left\{ 
\begin{array}{lr}
-\Delta {u_{0}(v)}=f, & x\in \Omega , \\ 
\partial _{\nu }u_{0}+\mathcal{A}_{1}\dfrac{a(x)}{a(x)+C_{0}}u_{0}(v)=v, & 
x\in \Gamma _{0}, \\ 
u_{0}(v)=0, & x\in \Gamma _{1}.%
\end{array}%
\right.
\end{equation*}
\end{theorem}

\noindent \textit{Proof.} Let us find the homogenized boundary value problem
for the function $u_{0}$. We set $\phi =W_{\varepsilon }\psi $ in the
integral identity for $u_{\varepsilon }$ and get 
\begin{equation*}
\lim_{{\varepsilon }\rightarrow 0}{\varepsilon }^{-\gamma
}\int_{S_{\varepsilon }}a(x)u_{\varepsilon }\phi {ds}=\mathcal{A}%
_{1}\int_{\Gamma _{0}}\frac{a(\hat{x})}{a(\hat{x})+C_{n}}u_{0}\phi {d\hat{x}}%
.\newline
\end{equation*}%
Hence, $u_{0}\in H^{1}(\Omega ,\Gamma _{1})$ satisfies integral identity 
\begin{equation}
\displaystyle\int_{\Omega }\nabla {u_{0}}\nabla \phi {dx}+\mathcal{A}%
_{1}\int_{\Gamma _{0}}\frac{a(\hat{x})}{a(\hat{x})+C_{n}}u_{0}\phi {d\hat{x}}%
=\int_{\Omega }f\phi {dx}-\frac{\eta }{N}\int_{\Gamma _{0}}P_{0}\phi {d\hat{x%
}},\newline
\label{Identity limit 32}
\end{equation}

\noindent for any $\phi \in H^{1}(\Omega ,\Gamma _{1})$. From here, we
derive that $u_{0}$ is a weak solution of the problem

\textrm{\textit{%
\begin{equation*}
\left\{ 
\begin{array}{lr}
-\Delta {u_{0}}=f, & x\in \Omega , \\ 
\partial _{\nu }u_{0}+\mathcal{A}_{n}\dfrac{a(x)}{a(x)+C_{n}}u_{0}=-\frac{%
\eta }{N}P_{0}, & x\in \Gamma _{0}, \\ 
u_{0}=0, & x\in \Gamma _{1}.%
\end{array}%
\right.
\end{equation*}%
}}

\noindent Next, let us obtain the limit problem satisfied by $P_{0}$. Let $%
\psi \in C^{\infty }(\overline{\Omega },\Gamma _{1})$. We take $\phi
=W_{\varepsilon }\psi $ as a test function in the integral identity for the
function $P_{\varepsilon }$ and get

\textrm{\textit{%
\begin{equation}
\int_{\Omega _{\varepsilon }}\nabla P_{\varepsilon }\nabla (W_{\varepsilon
}\psi )dx+\varepsilon ^{-\gamma }\int_{S_{\varepsilon }}a(x)P_{\varepsilon
}\psi ds=\int_{\Omega _{\varepsilon }}B(x)\nabla (u_{\varepsilon
}-u_{T})\nabla (W_{\varepsilon }\psi )dx.\   \label{Ident Peps 33}
\end{equation}%
}}

\noindent In order to pass to the limit, as $\varepsilon \rightarrow 0,$ in
the right-hand side of the identity (\ref{Identity limit 32}), we set $\phi
=u_{\varepsilon }\psi $ in the integral identity (\ref{Extension W ep 20})
and $\phi =h_{\varepsilon }\psi $ in the integral identity for $%
u_{\varepsilon }$. Then, subtracting one from the other, we get

\begin{equation}
\begin{array}{c}
\displaystyle \int _{\Omega _{\varepsilon }}B(x)\nabla {W_{\varepsilon }}%
\nabla (u_{\varepsilon }\psi )dx=\int\limits_{\Omega _{\varepsilon }}\nabla {%
h_{\varepsilon }}\nabla (u_{\varepsilon }\psi )dx-\int\limits_{\Omega
_{\varepsilon }}\nabla {u_{\varepsilon }}\nabla (h_{\varepsilon }\psi )dx+
\\ 
\displaystyle +\int\limits_{\Omega _{\varepsilon }}f{h_{\varepsilon }}\psi {%
dx}-\frac{\eta }{N}\int\limits_{\Gamma _{0}}P_{\varepsilon }h_{\varepsilon
}\psi {dx}.\newline
\end{array}
\label{ident P e limit 34}
\end{equation}

\noindent Identity (\ref{ident P e limit 34}) implies 
\begin{equation}
\begin{array}{c}
\displaystyle \lim _{\varepsilon \rightarrow 0}\int_{\Omega _{\varepsilon
}}B\nabla W_{\varepsilon }\nabla (u_{\varepsilon }\psi )dx= \\ 
\displaystyle =\int_{\Omega }\nabla {h_{0}}\nabla \psi {u_{0}}{dx}%
-\int_{\Omega }h_{0}\nabla {u_{0}}\nabla \psi {dx}+\int_{\Omega }f{h_{0}}%
\psi {dx}-\frac{\eta }{N}\int_{\Gamma _{0}}P_{0}h_{0}\psi d{\hat{x}}.%
\end{array}
\label{Identity B 35}
\end{equation}

\noindent We now use the fact that we already know the problems satisfied by 
$u_{0}$ and $h_{0}$. So, we get

\begin{equation}
\begin{array}{c}
\displaystyle\int_{\Omega }(\nabla {h_{0}},\nabla \psi )u_{0}dx-\int_{\Omega
}h_{0}\nabla {u_{0}}\nabla \psi {dx}=\int_{\Omega }\nabla {h_{0}}\nabla
(u_{0}\psi )dx-\int_{\Omega }\nabla {u_{0}}\nabla (h_{0}\psi )dx= \\[0.4cm]
\displaystyle-\frac{\mathcal{A}_{1}}{n}\int_{\Gamma _{0}}\frac{tr{B(\hat{x})}%
a(\hat{x})}{a(\widehat{x})+C_{n}}{u_{0}}\psi {d\hat{x}}-\int_{\Omega }f{h_{0}%
}\psi {dx}+\frac{\eta }{N}\int_{\Gamma _{0}}P_{0}\psi {d\hat{x}}%
\end{array}
\label{identity h0 and u 36}
\end{equation}

\noindent Comparing expressions (\ref{Identity B 35}) and (\ref{identity h0
and u 36}), we conclude

\textrm{\textit{%
\begin{equation}
\lim_{\varepsilon \rightarrow 0}\int_{\Omega _{\varepsilon }}B\nabla
W_{\varepsilon }\nabla (u_{\varepsilon }\psi )dx=-\frac{\mathcal{A}_{1}}{n}%
\int_{\Gamma _{0}}\frac{tr{B(\hat{x})}a(\hat{x})}{a(\hat{x})+C_{n}}u_{0}\psi
d\hat{x}.  \label{Limit B ans u0 37}
\end{equation}%
}}

\noindent We take $W_{\varepsilon }\psi $, with $\psi \in C^{\infty }(%
\overline{\Omega },\Gamma _{1})$, as a test function in the integral
identity for $P_{\varepsilon }$ and get

\begin{equation}
\int_{\Omega _{\varepsilon }}\nabla {W_{\varepsilon }}\nabla (P_{\varepsilon
}\psi )dx+{\varepsilon }^{-\gamma }\int _{S_{\varepsilon
}}a(x)P_{\varepsilon }\psi (x)ds=\int _{\Omega _{\varepsilon }}B\nabla {%
W_{\varepsilon }}\nabla (u_{\varepsilon }\psi )ds+\kappa _{\varepsilon },
\label{Limit W epsP 38}
\end{equation}%
where $\kappa _{\varepsilon }\rightarrow 0$ as $\varepsilon \rightarrow 0$.

\noindent From the definition of the function $W_{\varepsilon }$ and its
properties, we can transform the left-hand side of equality (\ref{Limit W
epsP 38}) in the following way 
\begin{equation*}
\sum_{j\in \Upsilon _{\varepsilon }}\int_{\partial {\ G_{\varepsilon }^{j}}%
\cup \partial {T_{{\varepsilon }/4}^{j}}}\partial _{\nu }w_{\varepsilon
}^{j}P_{\varepsilon }\psi {ds}+{\varepsilon }^{-\gamma }\int
_{S_{\varepsilon }}a(x)P_{\varepsilon }\psi {ds}=-\frac{\mathcal{A}_{1}}{n}%
\int _{\Gamma _{0}}\frac{tr{B(\hat{x})}{a(\hat{x})}}{a(\hat{x})+C_{n}}%
u_{0}\psi {d\hat{x}}+\kappa _{\varepsilon },\newline
\end{equation*}
where $\kappa _{\varepsilon }\rightarrow 0$ as $\varepsilon \rightarrow 0$.

\noindent The left-hand side of the last equality takes the form

\begin{equation}
{\varepsilon }^{-\gamma }\int_{S_{\varepsilon }}(a(x)+C_{n})P_{\varepsilon
}\psi {ds}+\sum _{j\in \Upsilon _{\varepsilon }}\int_{\partial {T_{{%
\varepsilon }/4}^{j}}}\partial _{\nu }w_{\varepsilon }^{j}P_{\varepsilon
}\psi {ds}+\theta _{\varepsilon },\newline
\label{Limit boundaries 39}
\end{equation}%
where $\theta _{\varepsilon }\rightarrow 0$ as $\varepsilon \rightarrow 0$.

\noindent From (\ref{Limit B ans u0 37})-(\ref{Limit boundaries 39}), we
deduce

\begin{equation}
\begin{array}{c}
\displaystyle \lim _{\varepsilon \rightarrow 0}{\varepsilon }^{-\gamma
}\int_{S_{\varepsilon }}a(x)P_{\varepsilon }\varphi {ds}=\mathcal{A}%
_{1}\int_{\Gamma _{0}}\frac{a(\hat{x})}{a(\hat{x})+C_{n}}P_{0}\varphi {d\hat{%
x}}-\frac{\mathcal{A}_{1}}{n}\int_{\Gamma _{0}}trB(\hat{x})\Bigl(\frac{a(%
\hat{x})}{a(\hat{x})+C_{n}}\Bigr)^{2}u_{0}\varphi {d\hat{x}},%
\end{array}
\label{Limit conclus 40}
\end{equation}

\noindent where $\varphi $ is an arbitrary function from $C^{\infty }(%
\overline{\Omega },\Gamma _{1})$. Passing to the limit in the integral
identity for $P_{\varepsilon }$ we get the theorem's statement.

\noindent Now we will find the limit as $\varepsilon \rightarrow 0$ of the
cost functional\textrm{\textit{%
\begin{equation*}
J_{\varepsilon }(v_{\varepsilon })=\frac{\eta }{2}\int_{\Omega _{\varepsilon
}}B\nabla (u_{\varepsilon }(v_{\varepsilon })-u_{T})\nabla (u_{\varepsilon
}(v_{\varepsilon })-u_{T})dx+\frac{N}{2}\Vert v_{\varepsilon }\Vert
_{L^{2}(\Gamma _{0})}^{2}.\newline
\end{equation*}%
}}

\noindent We take $\phi =u_{\varepsilon }$ in the integral identity for $%
P_{\varepsilon }$ and $\phi =P_{\varepsilon }$ in the integral identity for $%
u_{\varepsilon }$. Then, we subtract the one from the other and pass to the
limit as $\varepsilon \rightarrow 0$:

\textrm{\textit{%
\begin{equation}
A\equiv \lim_{{\varepsilon }\rightarrow 0}\int _{\Omega _{\varepsilon
}}B(x)\nabla (u_{\varepsilon }-u_{T})\nabla u_{\varepsilon }dx=\int_{\Omega
}fP_{0}dx-\frac{\eta }{N}\int_{\Gamma _{0}}P_{0}^{2}d\hat{x}.\newline
\label{Limit cost 42}
\end{equation}%
}}

\noindent Taking into account the integral identity for the limit of the
adjoint problem, we get

\begin{equation}
\begin{array}{c}
\displaystyle\mathcal{A}=\int_{\Omega }\nabla {u_{0}}\nabla {P_{0}}dx+%
\mathcal{A}_{1}\int_{\Gamma _{0}}\frac{a(\widehat{x})}{a(\widehat{x})+C_{n}}%
u_{0}P_{0}dx= \\[0.4cm]
\displaystyle=\int_{\Omega }B\nabla (u_{0}-u_{T})\nabla {u_{0}}dx+\frac{%
\mathcal{A}_{1}}{n}\int_{\Gamma _{0}}trB(\hat{x})\Bigl(\frac{a(\hat{x})}{a(%
\hat{x})+C_{n}}\Bigr)^{2}u_{0}^{2}d\hat{x}.\newline
\end{array}
\label{Cost limit 43}
\end{equation}

\noindent From here, we obtain (\ref{definition J0 41}). Since the trace is
a continuous operator on $H^{1}(\Omega _{\varepsilon },\Gamma _{1})$ we get
that $v_{\varepsilon }\rightharpoonup $ $v_{0}$ weakly in $L^{2}(\Gamma
_{0}) $. Since we have that $v_{0}=-\dfrac{\eta }{N}P_{0}$ on $\Gamma _{0}$
and this is the optimality condition associated to the control problem (\ref%
{Limit control problem}) (the proof is an easy variation of Proposition \ref%
{Propo1}), then $v_{0}$ is the unique optimal control problem associated to
the convex cost functional $J_{0}(v)._{\blacksquare }$

\bigskip

\begin{remark}
It is not too difficult to prove that the coupled system (\ref{Coupled limit
system 31}) has only one weak solution $(u_{0},P_{0})$. For instance, one
indirect proof can be obtained through the strict convexity of the
functional $J_{0}$. In particular, this implies that the weak convergence
obtained in Theorem \ref{Theorem 1} holds for any subsequences of the
original ones (since the limit $(u_{0},P_{0})$ is unique).
\end{remark}

\begin{remark}
Such as it is detailed explained in the book \cite{DGCSh 1 Book}, it can be
proved that the choice of the scales and parameters $\alpha =\gamma =\frac{%
n-1}{n-2}$ is the reason to get an anomalous boundary behavior on $\Gamma
_{0}$ (for bigger size of the elements of the lattice we get different
coefficients). This phenomenon was called in the literature as the
appearance of a \textquotedblleft strange term\textquotedblright\ and in
several papers it was associated to a certain \textquotedblleft
measure\textquotedblright\ $\mu .$ One of the merits of Theorem \ref{Theorem
Main} is to show that the \textquotedblleft strange term\textquotedblright\
is a certain completely identified function on $\Gamma _{0}.$
\end{remark}

\begin{remark}
Similar problems arise in many different applications, especially in the
field of Chemical Engineering (see, e.g., \cite{Gomez-Perez-sanchez} and
Chapter 5 of \cite{DGCSh 1 Book}). Some models in Climatology also use the
identification as a final boundary condition the limit of a thin layer on
which there are some suitable balances of differential equations and
transmissions conditions (see the so called energy balance models coupled
with a deep ocean in \cite{D-Hidal-Tello}). Problems quite similar to the
one considered in this paper arise also in Elasticity (see, e.g., \cite%
{Brezis-Caff-Friedman}).
\end{remark}

\begin{remark}
Many generalizations and applications seem possible and some of them will be
developed in some future works by the authors: i) Non-symmetrical shapes can
also be treated thanks to the techniques presented in \cite{DPSH Dokladi 5},
ii) The case of non-periodic lattices (under the assumption of
\textquotedblleft stationary and ergodic\textquotedblright\ random media)
can be considered as in the framwork traeted in \cite{Blanc-Lions}, \cite%
{Jikov-Kozlov-Oleinik}, \cite{Caffarelli Mellet} and \cite{Kruslov} (see
many other references in Appendice C of the book \cite{DGCSh 1 Book}). iii)
Optimal control problems for semilinear equations and/or nonlinear boundary
conditions could be approached as, for instance, in \cite{Conca et al}. iv)
The extension of the techniques of this paper can be also applied to the
consideration of several parabolic problems (see, e.g., Appendix A of \cite%
{DGCSh 1 Book}) and its references. v) by passing to the limit when
parameter $N\rightarrow +\infty $ it is possible to get some results on the
approximate controllability with internal observation (the $\ H^{1}$ norm of 
$(u_{\varepsilon }(v)-u_{T})$ can be made as small as wanted): see \cite%
{Glowinski he Lions}.
\end{remark}

\section{Convergence of the energy for the problem without control}

In this last Section, we consider the boundary value problem without any
control (i.e. problem (\ref{Eq 1}) with $v\equiv 0$) 
\begin{equation}
\left\{ 
\begin{array}{lr}
-\Delta {u_{\varepsilon }(v)}=f, & x\in \Omega _{\varepsilon }, \\ 
\partial _{\nu }u_{\varepsilon }(v)+{\varepsilon }^{-\gamma
}a(x)u_{\varepsilon }(v)=0, & x\in S_{\varepsilon }, \\ 
\partial _{\nu }u_{\varepsilon }(v)=0, & x\in \Gamma _{0}, \\ 
u_{\varepsilon }(v)=0, & x\in \Gamma _{1},%
\end{array}%
\right.  \label{Probl without}
\end{equation}%
where $f\in L^{2}(\Omega )$, $a(x)\in C^{\infty }(\overline{\Omega })$, $%
a(x)\geq a_{0}=const>0$. The homogenization techniques of previous sections
can be easily adapted to prove that $\widetilde{u_{\varepsilon }}%
\rightharpoonup u_{0}\,$weakly in $H^{1}(\Omega _{\varepsilon },\Gamma _{1})$
and that $u_{0}$ is a weak solution of the problem 
\begin{equation*}
\left\{ 
\begin{array}{lr}
-\Delta {u_{0}(v)}=f, & x\in \Omega , \\ 
\partial _{\nu }u_{0}+\mathcal{A}_{1}\dfrac{a(x)}{a(x)+C_{0}}u_{0}(v)=0, & 
x\in \Gamma _{0}, \\ 
u_{0}(v)=0, & x\in \Gamma _{1}.%
\end{array}%
\right.
\end{equation*}

Our main goal now is to prove that the consideration of an artificial
complementary system (formally corresponding to the case $v\equiv 0$ and $%
u_{T}\equiv 0$ and $B=I$) allows to prove the convergence of the
corresponding energies.

\begin{theorem}
\label{Them 2}Let $u_{\varepsilon }$ be the solution of (\ref{Probl without}%
) with $v\equiv 0$. Let $u_{0}\in H_{0}^{1}(\Omega )$ be the weak limit of
the extension $P_{\varepsilon }u_{\varepsilon }$. Then, we have the
convergence of the energy%
\begin{equation}
\int_{\Omega _{\varepsilon }}|\nabla u_{\varepsilon }|^{2}dx\rightarrow
\int_{\Omega }|\nabla u_{0}|^{2}dx+\mathcal{A}_{1}\int_{\Gamma _{0}}\Bigl(%
\frac{a(\hat{x})}{a(\hat{x})+C_{n}}\Bigr)^{2}u_{0}^{2}d\hat{x}.
\label{energy convergence}
\end{equation}
\end{theorem}

\noindent \textit{Proof}\textbf{. }We consider the auxiliary problem%
\begin{equation}
\left\{ 
\begin{array}{lr}
\Delta {P_{\varepsilon }}=\Delta u_{\varepsilon } & x\in \Omega
_{\varepsilon }, \\ 
\partial _{\nu }P_{\varepsilon }-\partial _{\nu }u_{\varepsilon }+{%
\varepsilon }^{-\gamma }a(x)P_{\varepsilon }=0, & x\in S_{\varepsilon }, \\ 
\partial _{\nu }P_{\varepsilon }-\partial _{\nu }u_{\varepsilon }=0, & x\in
\Gamma _{0}, \\ 
P_{\varepsilon }=0, & x\in \Gamma _{1}.%
\end{array}%
\right.
\end{equation}

\noindent As in the proof of Theorem \ref{Theorem Main}, we get that $%
\widetilde{P_{\varepsilon }}\rightharpoonup P_{0}\,$weakly in $H^{1}(\Omega
,\Gamma _{1})$ as $\varepsilon \rightarrow 0$, with $P_{0}$ the weak
solution of the problem%
\begin{equation}
\left\{ 
\begin{array}{lr}
\Delta {P_{0}}=\Delta u_{0} & x\in \Omega , \\ 
\partial _{\nu }P_{0}-\partial _{\nu }u_{0}+\mathcal{A}_{1}\dfrac{a(x)}{%
a(x)+C_{n}}P_{0}-\mathcal{A}_{1}\dfrac{a^{2}(x)}{(a(x)+C_{n})^{2}}u_{0}=0, & 
x\in \Gamma _{0}, \\ 
P_{0}=0, & x\in \Gamma _{1},%
\end{array}%
\right.  \label{Problem P without}
\end{equation}%
where $A_{1}(n)=(n-2)C_{0}^{n-2}\omega _{n}$, $C_{n}=\dfrac{n-2}{C_{0}}$ and 
$\omega _{n}$ is the surface area of the unit sphere in $\mathbb{R}^{n}$.

\noindent From the variational formulation of the problem (\ref{Probl
without}), taking ${P_{\varepsilon }}$ as a test function, we have

\begin{equation*}
\int_{\Omega _{\varepsilon }}\nabla {u_{\varepsilon }}\nabla P{_{\varepsilon
}}dx+{\varepsilon }^{-\gamma }\int_{S_{\varepsilon }}a(x)u_{\varepsilon
}P_{\varepsilon }ds=\int_{\Omega _{\varepsilon }}fP_{\varepsilon }dx.
\end{equation*}%
Similarly, from the variational formulation of the problem (\ref{Problem P
without}) on $P_{\varepsilon }$, taking ${u_{\varepsilon }}$ as a test
function, we derive 
\begin{equation*}
\int _{\Omega _{\varepsilon }}\nabla P{_{\varepsilon }}\nabla {%
u_{\varepsilon }}dx+{\varepsilon }^{-k}\int_{S_{\varepsilon
}}a(x)P_{\varepsilon }u_{\varepsilon }ds=\int_{\Omega _{\varepsilon
}}|\nabla {u_{\varepsilon }}|^{2}dx.
\end{equation*}%
Thus, we have 
\begin{equation*}
\begin{array}{l}
\displaystyle \int_{\Omega _{\varepsilon }}|\nabla {u_{\varepsilon }}%
|^{2}dx=\int_{\Omega _{\varepsilon }}fP{_{\varepsilon }}dx\rightarrow
\int_{\Omega }fP{_{0}}dx= \\[.5cm] 
\displaystyle =\int_{\Omega }\nabla {u_{0}}\nabla P{_{0}}dx+\mathcal{A}%
_{1}\int_{\Gamma _{0}}\frac{a(\hat{x})}{a(\hat{x})+C_{n}}u_{0}P_{0}d\widehat{%
x}= \\[.5cm] 
\displaystyle=\int_{\Omega }|\nabla {u_{0}}|^{2}dx+\mathcal{A}_{1}\int_{
\Gamma _{0}}\Bigl(\frac{a(\hat{x})}{a(\hat{x})+C_{n}}\Bigr)^{2}u_{0}^{2}d%
\hat{x},\newline
\end{array}%
\end{equation*}%
which ends the proof. $_{\blacksquare }$

\bigskip

\begin{remark}
It can be proved (for instance, by adapting the arguments presented in
Section 4.7.1.4 of \cite{DGCSh 1 Book}) that if we know that $u_{0}\in
W^{1,\infty }(\Omega )$ then we can get some results implying the strong
convergence of $u_{\varepsilon }$ plus a suitable \textquotedblleft
correction term\textquotedblright . Notice that the conclusion presented in
the proof of Theorem \ref{Them 2} follows different ideas.
\end{remark}

\bigskip

\textbf{Ackowledgements.} The research of J.I. D\'{\i}az was partially
supported the projects MTM2017-85449-P and PID2020-112517GB-I00 of the
DGISPI, Spain and the Research Group MOMAT (Ref. 910480) of the UCM.

\end{document}